\documentclass[12pt,letterpaper]{amsart}
\usepackage{euler, epic,eepic,latexsym, amssymb, amscd, amsfonts, xypic, url, color, epsfig}
\input xy
\xyoption{all}

\usepackage[OT2,T1]{fontenc}
\DeclareSymbolFont{cyrletters}{OT2}{wncyr}{m}{n}
\DeclareMathSymbol{\Sha}{\mathalpha}{cyrletters}{"58}
 
 \newlength{\baseunit}               
 \newcount{\numlines}                
\newcommand{\tpoint}[1]{\vspace{3mm}\par \noindent \refstepcounter{subsection}{\bf \thesubsection.} 
  {\bf #1. ---} }
\newcommand{\epoint}[1]{\vspace{3mm}\par \noindent \refstepcounter{subsection}{\thesubsection.} 
  {\bf #1.} }
\newcommand{\bpoint}[1]{\vspace{3mm}\par \noindent \refstepcounter{subsection}{\bf \thesubsection.} 
  {\bf #1.} }

\newcommand{\et}{\textrm{\'et}}
\newcommand{\rH}{{\rm H}}

\newcommand{\Z}{\mathbb{Z}}

\newcommand{\Q}{\mathbb{Q}}

\newcommand{\G}{\mathbb{G}}

\newcommand{\proj}{\mathbb P}

\newcommand{\Gal}{\operatorname{Gal}}

\newcommand{\Hom}{\operatorname{Hom}}

\newcommand{\Spec}{\operatorname{Spec}}

\newcommand{\GL}{\operatorname{GL}}
\newcommand{\Sym}{\operatorname{Sym}}

\newcommand{\Mat}{\operatorname{Mat}}
\newcommand{\Et}{\operatorname{Et}}
\newcommand{\Res}{\operatorname{Res}}
\newcommand{\Ker}{\operatorname{Ker}}
\newcommand{\Br}{\operatorname{Br}}
\newcommand{\Image}{\operatorname{Image}}

\newcommand{\hidden}[1]{\footnote{Hidden:  #1}}
\renewcommand{\hidden}[1]{}

\begin{document}
\pagestyle{plain}
\title{Splitting varieties for triple Massey products}

\author{Michael J. Hopkins}\author{Kirsten G. Wickelgren}\thanks{The first author is supported in part by the National Science Foundation under DMS-0906194 and by DARPA under HR0011-10-1-0054-DOD35CAP. The second author is supported by an American Institute of Mathematics five year fellowship.}
\address{Dept. of Mathematics, Harvard University, Cambridge~MA}
\email{mjh@math.harvard.edu}
\email{wickelgren@post.harvard.edu}
\date{October 2012.}
\begin{abstract}
We construct splitting varieties for triple Massey products. For $a,b,c \in F^*$ the triple Massey product $\langle a,b,c \rangle$ of the corresponding elements of $\rH^1(F, \mu_2)$ contains $0$ if and only if there is $x \in F^*$ and $y \in F[\sqrt{a}, \sqrt{c}]^*$ such that $b x^2 = N_{F[\sqrt{a}, \sqrt{c}]/F}(y)$, where $N_{F[\sqrt{a}, \sqrt{c}]/F}$ denotes the norm, and $F$ is a field of characteristic different from $2$. These varieties satisfy the Hasse principle by a result of D.B. Leep and A.R. Wadsworth. This shows that triple Massey products for global fields of characteristic different from $2$ always contain $0$.
\end{abstract}
\maketitle

{\parskip=12pt 

\section{Introduction}

Massey products measure information contained in a differential graded algebra which is lost by passing to its homology ring. For instance, the singular cochains of the complement of the Borromean rings have non-trivial Massey products, whence the differential graded algebra of cochains contains more information than the cohomology ring, and the Borromean rings are not equivalent to three unconnected circles. There are number theoretic analogues of this example, giving non-trivial Massey products in Galois cohomology  \cite{Morishita1} \cite[2.1.17]{Vogel_thesis} \cite[3.2-3.4]{Gaertner_thesis}.

Let $\eta$ be a functorial assignment of a set of cohomology classes $\eta_F \subseteq \rH^*(\Gal(F^s/F),\Z/2)$ to fields $F$ over a base field, where  $\rH^*(\Gal(F^s/F),\Z/2)$ denotes continuous group cohomology and $F^s$ denotes a separable closure of $F$. A splitting variety for $\eta$ is a scheme $X$ over the base field which has $F$-points if and only if there is an element of $\eta_F$ which vanishes.

This paper constructs a splitting variety for triple Massey products of elements of $\rH^1(\Gal(F^s/F), \Z/2)$ when $F$ is a field of characteristic $\neq 2$. Let $\kappa: F^* \to \rH^1(\Gal(F^s/F), \Z/2)$ denote the Kummer map, given by applying $\rH^*(\Gal(F^s/F), -)$ to the short exact sequence $$1 \to \mu_2 \to \G_m \stackrel{z \mapsto z^2}{\to} \G_m \to 1,$$ and identifying $\mu_2$ with $\Z/2$. For $a$,$b$,$c$ in $F^*$, let $\langle \kappa(a), \kappa(b), \kappa(c) \rangle$ denote the triple Massey product, when it is defined -- see \ref{MPdefn}.  $\langle \kappa(a), \kappa(b), \kappa(c) \rangle$ is a subset of $\rH^2(\Gal(F^s/F), \Z/2)$ which is a coset of a certain ideal called the indeterminacy, described in \ref{indeterminacy_remark}. Let $X(a,b,c)$ be the closed subscheme of $\G_{m} \times \mathbb{A}^4$ determined by the equation $$b x^2 = (y_1^2- a y_2^2 + c y_3^2 - ac y_{4}^2 )^2 - c(2 y_1 y_{3} - 2 a y_{2} y_{4} )^2 $$ where $x$ and $(y_1, y_2, y_3, y_4)$ are the coordinates on $\G_m$ and $\mathbb{A}^4$ respectively. The polynomial on the right hand side is the norm of $y_1 + y_2 \sqrt{a}  +  y_3 \sqrt{c} +y_4  \sqrt{a} \sqrt{c} $.

\tpoint{Theorem}\label{IntroSVthm}{\em $X(a,b,c)$ has an $F$-point if and only if $\langle \kappa(a), \kappa(b), \kappa(c) \rangle$ is defined and contains $0$.} 

An $R$-point of $X(a,b,c)$ for $a,b,c$ in $R^*$ implies that $\langle \kappa(a), \kappa(b), \kappa(c) \rangle$ is defined and contains $0$, for $R$ a ring with $2$ invertible and $\kappa$ taking values in the \'etale cohomology group $\rH^1(\Spec R, \mu_2)$ -- see Corollary \ref{Rpoint_implies_TripleMassey0}.

Theorem \ref{IntroSVthm}{ (= Theorem \ref{SplittingVarietyThm} below) is a special case of a more general construction, which could produce splitting varieties for order-$n$ Massey products of elements of $\rH^1(\Gal(F^s/F),\Z/2)$. See Section \ref{sectionSV}.

\epoint{Remark} $X(a,b,c)$ is the twist $  \{ x \times y \in \G_m \times \Res \G_m : b x^2 = \operatorname{Norm}(y) \}$ of the group scheme $$T = \{ x \times y \in \G_m \times \Res \G_m : x^2 = \operatorname{Norm}(y) \},$$ where $\Res$ denotes the restriction of scalars from $R[\alpha, \gamma]/\langle\alpha^2 - a, \gamma^2 -c \rangle$ to $R$ discussed in Section \ref{sectionSV}.  $\Res \G_m$ is isomorphic to the open subset of $\mathbb{A}^4 = \Spec R[y_1, y_2, y_3, y_4]$ where \begin{align*}\operatorname{Norm}(y_1, y_2, y_3, y_4) = \prod_{i,j = 0 ,1 } (y_1 + (-1)^i y_2 \alpha + (-1)^j y_3 \gamma + (-1)^{i+j}y_{4} \alpha \gamma) \\= (y_1^2- a y_2^2 + c y_3^2 - ac y_{4}^2 )^2 - c(2 y_1 y_{3} - 2 a y_{2} y_{4} )^2 \end{align*} is invertible. $T$ is an extension $$1 \to \mu_2 \to T \to \Res \G_m \to 1$$ of the torus $\Res \G_m$ by $\mu_2$. The $F$-points of $X(a,b,c)$ control the vanishing of $\langle \kappa(a),\kappa(b), \kappa(c) \rangle$, so we have that $\langle \kappa(a),\kappa(b), \kappa(c) \rangle$ contains $0$ if and only if there is $x \in F^*$ and $y \in F[\alpha, \gamma]/\langle\alpha^2 - a, \gamma^2 -c \rangle$ such that $b x^2 = \operatorname{Norm}(y)$.  Let $L= F[\sqrt{a}, \sqrt{c}]$ be the field extension of $F$ obtained by adjoining square roots of $a$ and $c$ in $F^s$. The  previous condition is equivalent to the existence of $x \in F^*$ and $y \in L^*$ such that $b x^2 = N_{L/F}(y)$ where $N_{L/F}$ denotes the norm (Corollary \ref{Nimage=Normimage}). 

We use the computed $X(a,b,c)$ to study triple Massey products. For instance, the equation for $X(a,b,c)$ gives another proof of the fact that the symmetric Massey product $\langle \kappa(a), \kappa(b), \kappa(a) \rangle$ contains $0$ whenever it is defined (Corollary \ref{<a,b,a>=0}). More interestingly, we show the vanishing of all triple Massey products of elements of $\rH^1(\Spec F, \Z/2)$, which are defined, for $F$ a global field. This is done as follows: the Hasse norm theorem mod squares of David Leep and Adrian Wadsworth  \cite{LW_transfer} \cite{LW_Hasse} shows that $X(a,b,c)$ satisfies the Hasse principle (Theorem \ref{HasseX}), and it follows that triple Massey products vanish globally if and only if they vanish everywhere locally. Furthermore, an alternate proof that triple Massey products satisfy the local-global principle would give an alternate proof of \cite[Thm 4.5]{LW_transfer} for n=2 which says that the Hasse norm theorem mod squares holds for extensions $F(\sqrt{a}, \sqrt{c})/F$. See Remark \ref{local_global_Galois_remark} for why this does not follow directly from the local-global principle for $\rH^2(-, \Z/2)$. In other words, the Hasse norm theorem mod squares for extensions  $F(\sqrt{a}, \sqrt{c})/F$ can be interpreted as the local-global principle for triple Massey products of elements of $\rH^1$ with $\Z/2$-coefficients. The Hasse norm theorem mod squares does not hold for all extensions of global fields\hidden{see introduction to \cite{LW_transfer}}.

By direct computation in Galois cohomology, one sees that triple Massey products vanish in $\Z/2$-cohomology of local fields (Lemma \ref{localM=0}). This shows:

\tpoint{Theorem}\label{Thm_globalM=0_intro}{\em Let $F$ be a global field of characteristic $\neq 2$ and $a,b,c \in F^*$. The triple Massey product $\langle \kappa(a),\kappa(b), \kappa(c) \rangle$ contains $0$ whenever it is defined.}

We end with a question. The Milnor conjecture proved by Voevodsky \cite[Thm 2.2]{Morel} \cite{V_Milnor_conjecture} says that the cohomology ring of $\Spec F$ is $$\rH^*(\Spec F, \Z/2) \cong \oplus_{m} \otimes_{i=1}^m F^*/ \langle 2, x \otimes (1-x)  \rangle.$$ It was shown in \cite{WKyoto} that many Massey products of $x$ and $1-x$ vanish in $\rH^2(\Spec F, \Z_{\ell})$, see loc. cit. Corollary 3.14.  Theorem \ref{Thm_globalM=0_intro} shows the vanishing of many triple Massey products with $\Z/2$-coefficients. We raise the question:

\bpoint{Question}{ Is $\mathcal{C}^*(\Spec F, \Z/2)$ formal?}

Since posting this paper on the arXiv, J\'an Min\'a\v{c} and Nguyen Duy Tan realized that joint work of Wengeng Gao, Leep, Min\'a\v{c}, and Tara Smith \cite{GLMS} implies that for any $a,b,c$ such that $\langle \kappa(a), \kappa(b), \kappa ( c) \rangle$ is defined, $X(a,b,c)$ has an $F$-point over any field $F$ of characteristic different from $2$. Moreover, they found an explicit simple rational point on these $X(a,b,c)$. Suresh Venapally independently observed this result as well. This greatly generalizes Theorem \ref{Thm_globalM=0_intro}.

Jochen G\"artner can construct relevant Massey products which do not vanish. Let  $(l_1,l_2,l_3) = (313,457,521)$ and consider Massey products in the differential graded algebra $\mathcal{C}^*(\pi_1(\Spec R), \Z/2)$ of $\Z/2$-cochains in continuous group cohomology of $\pi_1(\Spec R)$, where $R= \Spec \Z[\frac{1}{l_1}, \frac{1}{l_2}, \frac{1}{l_3}, \frac{1}{2}]$ and $\pi_1$ denotes the \'etale fundamental group. G\"artner can show that the Massey product $\langle \kappa(l_1), \kappa(l_2), \kappa(l_3) \rangle$ does not contain zero, and it follows that $X(313,457,521)$ has no $\Z[\frac{1}{313}, \frac{1}{457}, \frac{1}{512}, \frac{1}{2}]$-points. See Example \ref{Gartner}.


{\bf Acknowledgements}  It is a pleasure to thank Bjorn Poonen and Burt Totaro for very useful discussions. Bjorn Poonen sent us \cite{LW_transfer} and work on Theorem \ref{HasseX}. Burt Totaro sent \cite{Berhuy_Favi}, and we thank them both. We also wish to thank Jochen G\"artner for Example \ref{Gartner} and interesting correspondence, J\'an Min\'a\v{c} for interesting correspondence, and Aravind Asok, Jason Starr and all of the participants of AIM workshop {\em Rational curves and $\mathbb{A}^1$-homotopy theory} for discussing this problem with the second author.

\section{Splitting variety}\label{sectionSV}

Let $\mathcal{C}^*$ be a differential graded associative algebra with product $\cup$, differential $\delta: \mathcal{C}^* \to \mathcal{C}^{*+1}$, and homology $\rH^* = \ker \delta / \Image \delta$.

\bpoint{Definition of order-$n$ Massey products of elements of $\rH^1$}\label{MPdefn}  Choose an integer $n \geq 2$. A set $a_{ij}$ of elements of $\mathcal{C}^1$ for $1 \leq i < j \leq n+1$ and $(i,j) \neq (1, n+1)$ such that  $$\delta a_{ij} = \sum_{k = i+1}^{j-1} a_{ik} \cup a_{kj}$$ is called a {\em defining system} for the order-$n$ Massey product of the cohomology classes, denoted $a_1, a_2, \ldots, a_n$ respectively, represented by $a_{12}, a_{23} , ..., a_{n, n+1}$. The order-$n$ Massey product of  $a_1, a_2, \ldots, a_n$ is {\em defined} if there exists a defining system. The {\em order-$n$ Massey product} $\langle a_1, a_2, \ldots, a_n \rangle$ of  $a_1, a_2, \ldots, a_n$ with respect to the defining system $a_{ij}$ is  $$\langle a_1, a_2, \ldots, a_n \rangle_{\{ a_{ij} \}} =  \sum_{k = 2}^{n} a_{1k} \cup a_{k,n+1},$$ and when no defining system is specified, $\langle a_1, a_2, \ldots, a_n \rangle $ denotes the subset of $\rH^2$ consisting of the order-$n$ Massey products of  $a_1, a_2, \ldots, a_n$ with respect to all defining systems. The order-$3$ Massey product will be called the {\em triple Massey product}.

We describe Massey products in a little more detail. The order-$2$ Massey product is the cup product $\cup$. The next simplest is the triple Massey product of elements of $\rH^1$, which can be described as follows: suppose that $a,b,c \in \rH^1$ are such that $ab = bc = 0$. We can choose $A,B,C$ in $\mathcal{C}^1$ representing $a,b,c $ respectively. Since $ab = 0$, we can choose $E_{ab}$ such that $\delta E_{ab} = AB$, and similarly we can choose $E_{bc}$ such that $\delta E_{bc} = BC$. Note that $\delta(E_{ab}C - A E_{bc}) = 0$, whence $ E_{ab}C - A E_{bc}$ represents an element of $\rH^2$. The set of all $E_{ab}C + A E_{bc}$ obtained in this manner is the triple Massey product $\langle a,b,c \rangle \subseteq \rH^2$.  

Massey products arise naturally when attempting to classify differential graded algebras with a given cohomology ring. For instance, suppose $\rH^* \cong  \Z/2[a,b,c]/\langle ab, bc  \rangle$, where $a,b,c \in \rH^1$ are in degree $1$. Then $ E_{ab}C + A E_{bc}$ represents an element of $\rH^2$. Note that $\rH^2$ is a $\Z/2$-vector space with basis $\{a^2, b^2, c^2, ac\}$. By adjusting the choice of $E_{ab}$ and $E_{bc}$,  one can arrange $E_{ab}C + A E_{bc}$ to be cohomologous to $0$ or $b^2$, but not both, and this dichotomy changes the isomorphism class of $\mathcal{C}^*$.

\epoint{Remark on indeterminacy}\label{indeterminacy_remark} Choose $a_1$, $a_2$, and $a_3$ in $\mathcal{C}^*$. Suppose that $a_{ij}$ for $1 \leq i < j \leq 4$ and $(i,j) \neq (1, 4)$ is a defining system for the triple Massey product of $a_1$, $a_2$, $a_3$. It is straightforward to check that $\langle a_1, a_2, a_3 \rangle $ is the subset of $\rH^2$ given by \begin{equation*}\label{indeterminacy_eq}\langle a_1, a_2, a_3 \rangle = \langle a_1, a_2, a_3 \rangle_{\{a_{ij} \}} + a_1 \rH^1 + \rH^1 a_3.\end{equation*}\hidden{

To see this, first note that one may restrict to defining systems $a_{ij}$ with $a_{i ,i+1}$ fixed: choose $d_1$ in  $\mathcal{C}^1$. Then $$a_{12}' = a_{12} + \delta d_1, a_{23}' = a_{23} , a_{34}' = a_{34}$$ $$a_{13}' = a_{13} +  d_1 \cup a_{23},  a_{24}' = a_{24} ,$$  is a defining system for the triple Massey product of $a_1$, $a_2$, $a_3$. \begin{align*} \langle a_1, a_2, a_3 \rangle_{\{a_{ij}' \}}& = (a_{13} +  d_1 \cup a_{23}) \cup a_{34} + (a_{12} + \delta d_1) \cup a_{24}\\& = \langle a_1, a_2, a_3 \rangle_{\{a_{ij} \}} + d_1 \cup a_{23} \cup a_{34} + \delta d_1 \cup a_{24} \\& = \langle a_1, a_2, a_3 \rangle_{\{a_{ij} \}} + \delta(d_1 \cup a_{24}), \end{align*} whence we may fix the cocycle representative for $a_{1}$. It is straight-forward to check that we may likewise fix the cocycle representatives for $a_{2}$ and $a_3$.  \hidden{

Choose $d_2$ in $\mathcal{C}^1$. Then $$a_{12}' = a_{12}, a_{23}' = a_{23} + \delta d_2, a_{34}' = a_{34}$$ $$a_{13}' = a_{13} + a_{12} \cup d_2,  a_{24}' = a_{24} + d_2 \cup a_{34},$$  is a defining system for the triple Massey product of $a_1$, $a_2$, $a_3$. \begin{align*} \langle a_1, a_2, a_3 \rangle_{\{a_{ij}' \}}& = (a_{13} + a_{12} \cup d_2) \cup a_{34} + a_{12} \cup ( a_{24} + d_2 \cup a_{34}) \\ & =  \langle a_1, a_2, a_3 \rangle_{\{a_{ij} \}} + a_{12} \cup d_2 \cup a_{34} + a_{12} \cup d_2 \cup a_{34} & =  \langle a_1, a_2, a_3 \rangle_{\{a_{ij} \}}, \end{align*} whence we may fix the cocycle representative for $a_{23}$.}

Now for defining systems $a_{ij}$ and $a_{ij}'$ with $a_{i, i+1} = a_{i,i+1}'$, we have $$ \langle a_1, a_2, a_3 \rangle_{\{a_{ij} \}} -  \langle a_1, a_2, a_3 \rangle_{\{a_{ij}' \}} = (a_{13} - a_{13}') \cup a_{34} + a_{12}(a_{24} - a_{24}').$$ Since $\delta (a_{13} - a_{13}') = \delta(a_{24} - a_{24}') =0 $, it follows that $$ (a_{13} - a_{13}') \cup a_{34} + a_{12}(a_{24} - a_{24}') \in a_1 \rH^1 + \rH^1 a_3,$$ showing \eqref{indeterminacy_eq}.} The ideal $a_1 \rH^1 + \rH^1 a_3 \subseteq \rH^2$ is called the {\em indeterminacy}.

\epoint{Notation} Let $a_{ij}: \Mat_n (\Z/2) \to \Z/2$ be the function taking a $n \times n$ matrix with coefficients in $\Z/2$ to its $(i,j)$-entry. Let $$U_n = \{U \in \Mat_n(\Z/2): a_{ii}(U) = 1, a_{ij}(U) = 0 \textrm{ for all } i> j \}$$ be the group of unipotent $n \times n$ matrices with coefficients in $\Z/2$. 

\epoint{A result of Dwyer}\label{MPGroup_cohomology} Let $G$ be a group or a profinite group and let $\mathcal{C}^*$ denote the differential graded algebra of $\Z/2$-cochains in continuous group cohomology -- see for instance \cite[p. 14-15]{coh_num_fields} for the definition of inhomogeneous chains in continuous group cohomology. For homomorphisms $$a_i: G \to \Z/2$$ in $\mathcal{C}^1$ for $i=1,2,3$, the triple Massey product $ \langle a_1, a_2, a_3 \rangle$ contains $0$ if and only if there exists a homomorphism $G \to U_4$ such that the composition $$\xymatrix{G \ar[r] & U_4 \ar[rrr]^{a_{12} \times a_{23} \times a_{34}} &&& (\Z/2)^3}$$ is $a_1 \times a_2 \times a_3: G \to (\Z/2)^3$ by \cite[Thm 2.4]{Dwyer}.

We now place this result in the context relevant for this paper, and describe the procedure used to construct splitting varieties.

If $\mathcal{C}^*$ is the differential graded algebra of cochains on a space or pro-space $S$, say with $\Z/2$ coefficients, elements of $\rH^1$ correspond to $\Z/2$-torsors, which are equivalent to homomorphisms $\pi_1(S) \to \Z/2$. Any function $\pi_1(S) \to \Z/2$ gives a $1$-cochain on $K(\pi_1(S),1)$ and therefore determines an element of $\mathcal{C}^1$ by pulling back along the natural map from a pro-space to its Postnikov tower.

A $U_n$-torsor $P$ over $S$ gives rise to a homomorphism $\phi_P: \pi_1(S) \to U_n$. Composing with $a_{ij}$ yields $a_{ij} \phi_P$ in $\mathcal{C}^1$ whose boundary is computed by $$a_{ij} \phi_P (\gamma_1 \gamma_2) -  a_{ij} \phi_P (\gamma_1 ) - a_{ij} \phi_P (\gamma_2) = \sum_{i < k < j} a_{ik} \phi_P (\gamma_1 ) a_{jk} \phi_P (\gamma_2) $$ for $\gamma_1$, $\gamma_2$ in $\pi_1(S)$. Thus $$\delta  a_{ij} \phi_P = \sum_{i < k < j} a_{ik}  \phi_P \cup a_{kj} \phi_P.$$ Thus a $U_3$-torsor over $S$ such that $a_{12} \phi_P = a$ and $a_{23} \phi_P = b$ provides $E_{ab}$ as above by setting $E_{ab} = a_{13} \phi_P$. Similarly, a $U_4$-torsor such that  \begin{equation}\label{torsor_condition} a_{12} \phi_P = a, a_{23} \phi_P = b, a_{34} \phi_P = c\end{equation} produces the equation $$ \delta a_{14} \phi_P = E_{ab}c + a E_{bc}$$ in $\mathcal{C}^2$ where $E_{ab} = a_{13} \phi_P$ and $E_{bc} = a_{24} \phi_P$, from which it follows that $\langle a,b,c \rangle$ contains $0$. 

This algebraic manipulation is reversible and shows that the existence of a $U_4$-torsor satisfying \eqref{torsor_condition} is equivalent to $\langle a,b,c \rangle$ being defined and containing $0$ in the $\Z/2$-cochains on $K(\pi_1(S),1)$. Since $$\rH^2(\pi_1(S), \Z/2) \to \rH^2(S, \Z/2)$$ is injective\hidden{pf: let $$F \to S \to B \pi_1(S) $$ be the fiber sequence obtained from the map $S \to B \pi_1(S)$. The LES in $\pi_*$ shows that $\pi_1(F) = 0$, whence $\rH^1(F, \Z/2) = 0$. By the Serre spectral sequence $ H^i(\pi_1(S), H^j(F, \Z/2))\Rightarrow \rH^{i+j}(S, \Z/2)$, we have that the kernel of $\rH^2(\pi_1(S), \Z/2) \to \rH^2(S, \Z/2)$ is the image of $\rH^0(\pi_1(S), \rH^1(F, \Z/2)) = 0$, giving injectivity.} and the indeterminacy of any triple Massey product is contained in its image (see \ref{indeterminacy_remark}), this is equivalent to $\langle a,b,c \rangle$ being defined and containing $0$ in the $\Z/2$-cochains on $S$. The analogous statement holds for order-$n$ Massey products of elements of $\rH^1$ and $U_{n+1}$-torsors -- this is \cite[Thm 2.4]{Dwyer}. In particular, $U_4$-torsors imply the vanishing of associated Massey products in $\mathcal{C}^*$, and more generally, unipotent representations of the fundamental group give information about the differential graded algebra of cochains.

To a scheme $S$, we can associate the \'etale topological type $\Et(S)$ and the ind-differential graded algebra of $\Z/2$-cochains on $\Et(S)$. The associated cohomology ring is $\rH^*(S, \Z/2)$, where $\rH^*$ denotes \'etale cohomology by \cite[Prop 5.9]{Friedlander}. By the above, to every $U_4$-torsor over $S$, there is an associated Massey product which vanishes. Note that in this context, $U_4$ is considered as a constant group scheme. Furthermore, for any triple $a$,$b$,$c$ of elements of $\rH^1(S, \Z/2)$, the vanishing of $\langle a,b,c \rangle$ is equivalent to the existence of a $U_4$-torsor whose push-forward along $$a_{12} \times a_{23} \times a_{34}: U_4 \to \Z/2 \times  \Z/2 \times \Z/2$$ is classified by $a \times b \times c$.  Thus a universal $U_4$-torsor with prescribed push-forward gives a splitting variety for a triple Massey product. More generally, a universal $U_{n+1}$-torsor with prescribed push-forward gives a splitting variety for an order-$n$ Massey product of elements of $\rH^1$.

In topological spaces, a universal $U_{n+1}$-torsor is produced by the quotient of a contractible space by a free action of $U_{n+1}$.  In schemes, there are $U_{n+1}$-torsors which are universal for $U_{n+1}$-torsors over $\Spec F$, where $F$ ranges over all field extensions of some base field, in the sense that a $U_{n+1}$-torsor $\mathcal{V} \to \mathcal{X}$ is universal when the map from $\mathcal{X}(F)$ to isomorphism classes of $U_{n+1}$-torsors over $\Spec F$ is surjective. To emphasize the difference between requiring surjectivity and bijectivity, one could call such a torsor versal, instead of universal, but we won't use this convention here. A universal torsor can be constructed from a linear action with trivial stabilizers on an open subscheme of affine space. This uses Hilbert 90 and an understanding of points over a field. See for instance \cite[Ex 5.4 p 12]{GMS} \cite[Prop 4.11]{Berhuy_Favi}.

A universal $U_{n+1}$-torsor with controlled push-forward along $$q = a_{12} \times a_{23} \times \ldots \times a_{n, n+1} $$ can be constructed from a universal $U_{n+1}$-torsor and a universal $(\Z/2)^{n}$-torsor. Let  $\mathcal{V} \to \mathcal{V} / U_{n+1}$ and $\mathcal{L} \to \mathcal{L}/ (\Z/2)^{n}$ be universal $U_{n+1}$ and $(\Z/2 )^{n}$-torsors respectively. Endowing $\mathcal{V} \times \mathcal{L}$ with the diagonal action of $U_{n+1}$, form $\mathcal{X} = (\mathcal{V} \times \mathcal{L})/ U_{n+1} $. The quotient $\mathcal{V} \times \mathcal{L} \to \mathcal{X} $ is then a universal $U_{n+1}$-torsor, and there is a tautological map of $(\Z/2)^{n}$-torsors $$ \xymatrix{q_* (\mathcal{V} \times \mathcal{L}) \ar[r] \ar[d] & \mathcal{L} \ar[d] \\ \mathcal{X} \ar[r] &\mathcal{L}/(\Z/2)^{n} }.$$ For a chosen $(\Z/2)^{n}$-torsor over $\Spec F$, the universality of $\mathcal{L} \to \mathcal{L}/ (\Z/2)^{n}$ gives an $F$-point $x$ of $\mathcal{L}/(\Z/2)^{n}$ such that pull-back along $x$ is the chosen torsor. The restriction of $\mathcal{V} \times \mathcal{L}$ to the fiber $\mathcal{X} (x)$ of $\mathcal{X}$ over $x$ is a universal $U_{n+1}$-torsor with prescribed push-forward. For $x$ such that the corresponding $(\Z/2)^{n}$-torsor is classified by $\kappa(a_1) \times \kappa(a_2) \times \ldots \times \kappa(a_n)$, the resulting $\mathcal{X} (x)$ has the property that $S$-points imply the vanishing of $\langle \kappa(a_1), \kappa(a_2), \ldots , \kappa(a_n) \rangle$ in $\rH^2(S, \Z/2)$, and this implication becomes an equivalence for $S = \Spec F$, for $F$ a field. The constructed torsor is not unique.

One can obtain defining equations for universal $G$-torsors (with or without controlled push-forward) for $G$ a finite $2$-group  as follows: if $R$ is a ring with $2 \in R^*$, then an $R$-module $A$ with an action of $(\Z/2)^m$ decomposes into a direct sum of simultaneous eigenspaces,\hidden{Let $M$ be the $R$-module. For $m=1$, let $tau$ generate $\Z/2$. Let $M^+ = \{m \in M: \tau m = m \}$ and $M^- = \{m \in M: \tau m = -m \}$. Then the natural map $M^+ \oplus M^{-} \to M$ is an isomorphism because it is injective as $M^+ \cap M^- = 0$ and surjective because $m = 1/2(m + \tau m) + 1/2(m - \tau m)$. For $m>1$, the decomposition for the first $\tau$ is respected by the action of $\Z/2^{m-1}$, so proceed by induction.} allowing us to compute the fixed elements $A^{(\Z/2)^m}$.  Filtering $G$ as $$1 \subset G_0 \subset G_1 \subset \ldots \subset G$$ such that $G_n/G_{n-1} \cong \Z/2^{m_n}$ and successively computing fixed elements of $A^{G_{n-1}}$ under the action of $G_n/G_{n-1} \cong \Z/2^{m_n}$, where $A$ is a ring of functions, produces the ring of functions of a universal torsor $\mathcal{X} =  (\mathcal{V} \times \mathcal{L})/G$.

\epoint{Example} For $G = U_3$, the Brauer-Severi variety $av^2 + b w^2 = u^2$ is obtained in this manner as a splitting variety for the cup-product $\langle \kappa(a), \kappa(b) \rangle = \kappa(a) \cup \kappa(b)$.  To see this:   consider the subgroup of $U_3$ of matrices with $a_{12} = 0$. The rank $1$ representation of this subgroup defined by a matrix $g$ acting by $(-1)^{a_{13}(g)}$ gives a rank $2$ representation of $U_3$ by coinduction. Let $\mathcal{V}$ be $\mathbb{A}^2 - \{ 0\}$ where $U_4$ acts via the coinduced representation. More explicitly, identify $\mathbb{A}^2$ with $\mathbb{A}^2= \Spec R[x,y] $. For $i < j$, let $E_{ij} \in U_3$ be the matrix such that $a_{ij}= 1$ and the other non-diagonal entries are $0$. Then $E_{13}$ acts by $E_{13} (x) = -x$ and $E_{13}(y) = -y$, $E_{12}$ acts by switching $x$ and $y$, and $E_{23}$ acts by $E_{23} (x) = -x$ and $E_{23}(y) = y$. Let $\mathcal{L}$ be $\G_m \times \G_m = \Spec R[\alpha, \alpha^{-1}, \beta, \beta^{-1}]$. Let $e_{12}$ and $e_{23}$ be a basis for $\Z/2 \times \Z/2$ and let $\Z/2 \times \Z/2$ act on $\mathcal{L}$ by $e_{12} (\alpha) = - \alpha$, $e_{12} (\beta) =  \beta$, $e_{23} (\alpha) =  \alpha$, and $e_{23} (\beta) = - \beta$. The quotient $(\mathbb{A}^2 \times \mathcal{L})/ \langle E_{14} \rangle$ is $\Sym^2 \mathbb{A}^2 \times \mathcal{L}$, which has ring of functions $R[x^2, y^2, xy,  \alpha, \alpha^{-1}, \beta, \beta^{-1}]$. The symmetric square $\Sym^2 \mathbb{A}^2$ with its action of $U_3/\langle E_{14} \rangle \cong (\Z/2) E_{12} \times (\Z/2) E_{23}$  has simultaneous eigenfunctions $x^2 - y^2$, $x^2 + y^2$, $ xy$ with eigenvalues $(-1,1)$, $(1,1)$, and $(1,-1)$ respectively. It follows that the ring of functions of $(\mathbb{A}^2 \times \mathcal{L})/ U_3$ is $R[(x^2 - y^2)/ \alpha, x^2 + y^2, xy/ \beta, \alpha^2, \alpha^{-2}, \beta^2, \beta^{-2}]$. The splitting variety  $X(a,b)$ is defined as the pullback of $(\mathbb{A}^2 \times \mathcal{L})/ U_3 \to \mathcal{L}/(\Z/2 \times \Z/2)$ by an $R$-point of $\mathcal{L}/(\Z/2 \times \Z/2)$ classifying $\kappa(a) \cup \kappa(b)$. The quotient $\mathcal{L}/(\Z/2 \times \Z/2)$ has function ring $R[\alpha^2, \alpha^{-2}, \beta^2, \beta^{-2}]$ and $\kappa(a) \cup \kappa(b)$ corresponds to the $R$-point obtained by letting $a = \alpha^2$ and $b = \beta^2$. We see that $v \mapsto (x^2 - y^2)/ \alpha$, $w \mapsto 2 xy /\beta$, and $u \mapsto x^2 + y^2$ defines an isomorphism from $X(a,b)$ to the Brauer-Severi variety. 

We now turn to the main theorem of this paper, which is the computation of a universal $U_4$-torsor. 

Let $R$ be a ring such that $2$ is invertible in $R$.

Let $\kappa: R^* \to \rH^1(\Spec R, \mu_2)$ denote the Kummer map obtained by applying  $\rH^*(\Spec R, -)$ to $$1 \to \mu_2 \to \G_m \stackrel{z \mapsto z^2}{\to} \G_m \to 1,$$ where $\rH^*$ denotes \'etale cohomology.

We will use multiple copies of $\G_{m,R}$. For notational convenience, let $L_{x} = \Spec R[x, x^{-1}] \cong \G_{m,R}$ denote a copy of $\G_{m,R}$ with distinguished function $x$. Let $$S = L_a \times L_c = \Spec R[a,c, a^{-1}, c^{-1}]$$ $$S' =L_{\alpha} \times L_{\gamma} = \Spec R[\alpha,\gamma, \alpha^{-1}, \gamma^{-1}]$$and $S' \to S$ be the degree $4$ finite \'etale cover $$a \mapsto \alpha^2, c \mapsto \gamma^2.$$

Let $\Res_{S'/S} \G_{m,S'}$ denote the restriction of scalars \cite[\S 7.6]{BLR}. Form the $L_a \times L_b \times L_c = \Spec R[a, b, c,  a^{-1}, b^{-1},c^{-1}]$-scheme $$\Res_{S'/S} \G_{m,S'} \times \G_{m, L_{b}}$$ obtained by taking the product over $ \Spec R$ of $\Res_{S'/S} \G_{m, S'}$ and $\G_{m,L_b}$. Let $x$ denote a function on $\G_{m, L_{b}}$ such that $\G_{m, L_b} = \Spec R[b,b^{-1},x,x^{-1}]$. Let $N_{S'/S}$  denote the norm on $\Res_{S'/S} \G_{m,S'}$ $$N_{S'/S}(y) = \prod_{i,j = 0 ,1 } (y_1 + (-1)^i y_{\alpha} \alpha + (-1)^j y_{\gamma} \gamma + (-1)^{i+j}y_{\alpha \gamma} \alpha \gamma),$$ where $$y = y_1 + y_{\alpha} \alpha + y_{\gamma} \gamma + y_{\alpha \gamma} \alpha \gamma .$$

\bpoint{Definition} Let $\mathcal{X} \subset \Res_{S'/S} \G_{m,S'} \times \G_{m, L_b}$ be the closed subscheme determined by  $$b x^2 = N_{S'/S}(y).$$ 

Let $\pi: \mathcal{X} \to L_a \times L_b \times L_c = \Spec R[a,c, b, a^{-1}, b^{-1},c^{-1}]$ denote the structure map. 

\bpoint{Definition} For $(a,b,c) \in (L_a \times L_b \times L_c)(R)$, let $X(a,b,c)$ denote the $R$-scheme obtained by pulling-back $\mathcal{X}$ along $(a,b,c)$. So $X(a,b,c)$ is the affine scheme with coordinate ring $$R[x,x^{-1}, y_1, y_{\alpha}, y_{\gamma}, y_{\alpha \gamma}]/ \langle b x^2 - N_{S'/S} (y_1 + y_{\alpha} \alpha + y_{\gamma} \gamma + y_{\alpha \gamma} \alpha \gamma) \rangle,$$ where we note that $$N_{S'/S} (y_1 + y_{\alpha} \alpha + y_{\gamma} \gamma + y_{\alpha \gamma} \alpha \gamma) = (y_1^2- a y_{\alpha}^2 + c y_{\gamma}^2 - ac y_{\alpha \gamma}^2 )^2 - c(2 y_1 y_{\gamma} - 2 a y_{\alpha} y_{\alpha \gamma} )^2$$ is indeed an element of $R[x,x^{-1}, y_1, y_{\alpha}, y_{\gamma}, y_{\alpha \gamma}]$. 

Let $q:U_4 \to (\Z/2)^3$ be the quotient group homomorphism $q = a_{12} \times a_{23} \times a_{34}$. Adopt the convention that a group also denotes the corresponding constant group scheme over $R$. 

For $P$ a $U_4$-torsor over a scheme $X$, the homomorphism $q$ determines a $(\Z/2)^3$-torsor over $X$  $$q_* P = P \times^{U_4} (\Z/2)^3 :=  P \times (\Z/2)^3 / U_4$$  or equivalently $q_*: \rH^1(-, U_4) \to \rH^1(-, (\Z/2)^3)$ is the map induced by $q$. 

Let $\mathcal{L} = L_{\alpha} \times L_{\beta} \times L_{\gamma}$, and let $\mathcal{L} \to L_a \times L_b \times L_c$ be the $(\Z/2)^3$-torsor given by $$a \mapsto \alpha^2, b \mapsto \beta^2, c \mapsto \gamma^2 .$$

Let $U_4$ act on $\mathcal{L}$ by $g^* \alpha = (-1)^{a_{12}(g)} \alpha$, $g^* \beta = (-1)^{a_{23}(g)} \beta$, and $g^* \gamma = (-1)^{a_{34}(g)} \gamma$ for $g$ in $U_4$.

For $i < j$, let $E_{ij}$ in $U_4$ denote the matrix such that $a_{ij}(E_{ij}) = 1$ and $a_{kl}(E_{ij}) = 0$ for $k \neq l$. 

Consider the following representation $V$ of $U_4$: let $\tilde{H}$ denote the subgroup of $U_4$ of matrices such that $a_{12} = 0$ and $a_{13} = 0$. Note that $\tilde{H}$ decomposes as the product $\tilde{H} = Z \times H$ of $Z = \{E_{14}, 1\}$ and $H= \Ker (a_{14}: \tilde{H} \to \Z/2)$. The notation is chosen to recall that $H$ is isomorphic to the Heisenberg group.

The homomorphism $a_{14}$ gives a $1$-dimensional representation $\sigma_{14}$ of $\tilde{H}$, where $h \in \tilde{H}$ acts on $R$ by multiplication by $(-1)^{a_{14}(h)}$. Let $V$ denote the coinduced  representation $$V = \operatorname{Ind}_{\tilde{H}}^{U_4} \sigma_{14} = R[U_4] \otimes_{R[\tilde{H}]} \sigma_{14}.$$

Let $\mathcal{V} \to \Spec R$ be the vector bundle associated to the representation $V$. Namely, let $V^* = \Hom_R (V, R)$ denote the linear dual of $V$ and let $\mathcal{V} =\Spec \Sym V^*$. We specify coordinates on $\mathcal{V}$ as follows. $V$ is a free $R$ module of rank $4$ with a basis corresponding to the cosets $\tilde{H}$, $E_{12} \tilde{H}$, $E_{13} \tilde{H}$, $E_{12}E_{13}\tilde{H}$ of $\tilde{H}$ in $U_4$. Let $\{ u_1, u_2, u_3, u_4 \}$ denote the basis of $V$ defined by $ u_1=\tilde{H} + E_{13} \tilde{H}$, $u_2 = \tilde{H} -E_{13} \tilde{H}$, $u_3 = E_{12} \tilde{H} + E_{12}E_{13} \tilde{H}$, and $u_4 = E_{12} \tilde{H} - E_{12}E_{13}  \tilde{H}$.  We use here that $2$ is invertible in $R$. For $i = 1,2,3,4$, let $u_i^*: V \to R$ denote the basis dual to $\{ u_1, u_2, u_3, u_4\}$, i.e. the functions such that $v = u_1^*(v) u_1 + u_2^*(v) u_2 + u_3^*(v) u_3 + u_4^*(v) u_4$ for all $v$ in $V$. Then $\mathcal{V} = \Spec R[u_1^*, u_2^*, u_3^*, u_4^*]$. Note that $U_4$ acts on $\mathcal{V}$ through the representation $V$. 

Let $(\mathcal{V} \times \mathcal{L})/U_4$ denote the quotient scheme of $\mathcal{V} \times \mathcal{L}$ by $U_4$ as in \cite[V Prop. 1.1]{sga1}. Other quotient schemes will be denoted similarly.

By a slight abuse of notation, let $N_{S'/S} (y_1 + y_{\alpha} \alpha + y_{\gamma} \gamma + y_{\alpha \gamma} \alpha \gamma) $ denote the element of $$R[a,b,c, y_1, y_{\alpha}, y_{\gamma}, y_{\alpha \gamma}] $$ given by setting $a = \alpha^2$, $b = \beta^2$, $c = \gamma^2$ and expanding $$ \prod_{i,j = 0 ,1 } (y_1 + (-1)^i y_{\alpha} \alpha + (-1)^j y_{\gamma} \gamma + (-1)^{i+j}y_{\alpha \gamma} \alpha \gamma) = (y_1^2- a y_{\alpha}^2 + c y_{\gamma}^2 - ac y_{\alpha \gamma}^2 )^2 - c(2 y_1 y_{\gamma} -2 a y_{\alpha} y_{\alpha \gamma} )^2.$$

\tpoint{Lemma}\label{coords_for_VxL/U_4}{\em There is an isomorphism $$ (\mathcal{V} \times \mathcal{L})/U_4 \cong \Spec R[a, a^{-1}, b, b^{-1}, c, c^{-1}, y_1, y_{\alpha}, y_{\gamma}, y_{\alpha \gamma}, x] / \langle b x^2 - N_{S'/S} (y_1 + y_{\alpha} \alpha + y_{\gamma} \gamma + y_{\alpha \gamma} \alpha \gamma) \rangle $$ sending $$x \mapsto 2^4 u_1^* u_2^* u_3^* u_4^*/(\beta) ,$$}

\begin{proof}
It is straightforward to verify that \begin{enumerate} \item \label{E13E14E24action} $\{ u_1, u_2, u_3, u_4 \}$ are simultaneous eigenvectors for $( E_{13}, E_{24},E_{14})$ with eigenvalues $(1,1, -1)$, $(-1,1,-1)$, $(1,-1,-1)$, $(-1,-1,-1)$ respectively.
\item \label{E12action} $E_{12}$ acts on $\{ u_1, u_2, u_3, u_4 \}$ by the permutation $(u_1,u_3)(u_2,u_4)$.
\item \label{E23action} $\{ u_1, u_2, u_3, u_4 \}$ are eigenvectors for $E_{23}$ with eigenvalues $1,1,1,-1$ respectively.
\item \label{E34action} $E_{34}$ acts on $\{ u_1, u_2, u_3, u_4 \}$ by the permutation $(u_1,u_2)(u_3,u_4)$.
\end{enumerate}

Since $\langle E_{13}, E_{24}, E_{14} \rangle$ acts trivially on $\mathcal{L}$, $$ (\mathcal{V} \times \mathcal{L})/ \langle E_{13}, E_{24}, E_{14} \rangle  = \Spec R[u_1^*, u_2^*, u_3^*, u_4^*]^{\langle E_{13}, E_{24}, E_{14} \rangle} \times \mathcal{L}.$$ $R[u_1^*, u_2^*, u_3^*, u_4^*]$ decomposes into simultaneous eigenspaces for $\langle E_{13}, E_{24}, E_{14} \rangle \cong (\Z/2)^3$.  The $(1,1,1)$-eigenspace is equal to the subalgebra $R[(u_1^*)^2, (u_2^*)^2, (u_3^*)^2, (u_4^*)^2, u_1^*u_2^*u_3^*u_4^* ]$ by \eqref{E13E14E24action}.

$R[u_1^*, u_2^*, u_3^*, u_4^*]^{ \langle E_{13}, E_{24}, E_{14} \rangle} = R[(u_1^*)^2, (u_2^*)^2, (u_3^*)^2, (u_4^*)^2, u_1^*u_2^*u_3^*u_4^*]$ decomposes into simultaneous eigenspaces for $U_4 / \langle E_{13}, E_{24}, E_{14} \rangle \cong (\Z/2)^3$. A $\Z/2$ basis for $U_4 / \langle E_{13}, E_{24}, E_{14} \rangle$ is $\{E_{12}, E_{23}, E_{34} \}$.

By \eqref{E12action}  \eqref{E23action}  \eqref{E34action}, $d_1$,$d_2$,$d_3$,$d_4$  defined by

$$d_1= (u_1^*)^2 + (u_2^*)^2 +( u_3^*)^2 + (u_4^*)^2$$
$$d_2= (u_1^*)^2 + (u_2^*)^2 - (u_3^*)^2 -( u_4^*)^2$$
$$d_3= (u_1^*)^2 - (u_2^*)^2 + (u_3^*)^2 - (u_4^*)^2$$
$$d_4= (u_1^*)^2 - (u_2^*)^2 - (u_3^*)^2 + (u_4^*)^2$$

\noindent are simultaneous eigenvectors for $\{E_{12}, E_{23}, E_{34} \}$ with eigenvalues $(1,1,1)$, $(-1,1,1)$, $(1,1,-1)$, $(-1,1,-1)$ respectively. Since $d_1$,$d_2$,$d_3$,$d_4$, and $ u_1^*u_2^*u_3^*u_4^*$ are algebra generators and simultaneous eigenvectors, they determine the eigenspace decomposition of $$R[u_1^*, u_2^*, u_3^*, u_4^*]^{ \langle E_{13}, E_{24}, E_{14} \rangle}.$$

The eigenspace decomposition of $R[\alpha, \alpha^{-1}, \beta, \beta^{-1}, \gamma, \gamma^{-1}]$ is apparent, giving the eigenspace decomposition of $$R[ (u_1^*)^2, (u_2^*)^2, (u_3^*)^2, (u_4^*)^2, u_1^*u_2^*u_3^*u_4^*, \alpha, \alpha^{-1}, \beta, \beta^{-1}, \gamma, \gamma^{-1}].$$

It follows that $$(\mathcal{V} \times \mathcal{L})/U_4 \cong \Spec R[\alpha^2, \alpha^{-2},\beta^2, \beta^{-2}, \gamma^2, \gamma^{-2}, d_1, d_2/\alpha, d_3/\gamma, d_4/ (\alpha \gamma), u_1^* u_2^* u_3^* u_4^*/(\beta)],$$ where the ring on the right hand side denotes the subalgebra of $$R[  u_1^*,u_2^*,u_3^*,u_4^*, \alpha, \alpha^{-1}, \beta, \beta^{-1}, \gamma, \gamma^{-1}]$$ generated by the listed elements.

It follows that sending $$x \mapsto 2^4 u_1^* u_2^* u_3^* u_4^*/(\beta) ,$$ $$y_1 \mapsto d_1,$$ $$y_{\alpha} \mapsto d_2/\alpha,$$ $$y_{\gamma} \mapsto d_3/\gamma,$$ $$y_{\alpha \gamma} \mapsto d_4 /(\alpha \gamma),$$ and $a \mapsto \alpha^2$, $b \mapsto \beta^2$, $c \mapsto \gamma^2$ defines the required isomorphism.

\end{proof}

\tpoint{Theorem}\label{have_torsors}{\em There is a $U_4$-torsor $\G_{m, R}^4 \times \mathcal{L} \to \mathcal{X}$ and an isomorphism of $(\Z/2)^3$-torsors $q_* (\G_{m, R}^4 \times \mathcal{L}) \to \pi^* \mathcal{L}$.}

\begin{proof}
Let $\G_{m,R}^4 \to \mathcal{V}$ be the open immersion onto the complement of the zero locus of $u_1^* u_2^* u_3^* u_4^* $. By Lemma \ref{coords_for_VxL/U_4} and \cite[V Cor1.4]{sga1}, $$(\G_{m,R}^4 \times \mathcal{L})/U_4 \cong \mathcal{X}.$$ 


The resulting map $f: (\G_{m,R}^4 \times \mathcal{L}) \to \mathcal{X}$ will be shown to be the claimed $U_4$-torsor.

Note that $f$ is finite-type. Thus by \cite[V Prop 1.1(i)]{sga1}, $f$ is finite.

Let $x$ be a point of $\G_{m,R}^4 \times \mathcal{L}$, and let $I(x)$ be the inertia group of $x$ i.e. the subgroup of $U_4$ of elements which stabilize $x$ and which act trivially on the residue field $k(x)$. Note that $\alpha$, $\beta$, and $\gamma$ are non-zero elements of $k(x)$. Since $g^* \alpha = (-1)^{a_{12}(g)} \alpha$ for $g$ in $U_4$, we have that $a_{12}(g) = 0$ for $g$ in $I(x)$. Similarly, $a_{23}(g) = a_{34}(g) = 0$.  The elements of $U_4$ such that $a_{12}(g) = a_{23}(g)= a_{34}(g) = 0$ form the subgroup $\langle E_{13}, E_{14}, E_{24} \rangle \cong (\Z/2)^3$. Since $u_1^*$ determines a non-zero element of $k(x)$ and for all $g$ in $\langle E_{13}, E_{14}, E_{24} \rangle$ the action of $g$ on $u_1^*$ is $g^* u_1^* = (-1)^{a_{14}(g)} u_1^*$ by \eqref{E13E14E24action}, we have that $a_{14}(g) = 0$ for $g$ in $I(x)$. Since $u_2^*$ determines a non-zero element of $k(x)$, we have by \eqref{E13E14E24action} that $a_{13}(g) + a_{14}(g) = 0$, whence $a_{13}(g) = 0$  for $g$ in $I(x)$.  Since $u_3^*$ determines a non-zero element of $k(x)$, we have by \eqref{E13E14E24action} that $a_{24}(g) + a_{14}(g) = 0$, whence $a_{24}(g) = 0$  for $g$ in $I(x)$. Thus $I(x)$ is $\{1\}$. 

Since the inertia groups at all points of $\G_{m,R}^4 \times \mathcal{L}$ are trivial and $f$ is finite, it follows that $f$ is a $U_4$-torsor by \cite[V Prop 2.6 (i) (ii)]{sga1}.

The map $\G_{m,R}^4 \times \mathcal{L} \times (\Z/2)^3 \to \mathcal{L}$ defined by projecting to $\mathcal{L} \times (\Z/2)^3$ and composing with the right multiplication $\mathcal{L} \times (\Z/2)^3 \to \mathcal{L}$ determines a well-defined map $$q_* f = (\G_{m,R}^4 \times \mathcal{L}) \times^{U_4} (\Z/2)^3  \to \mathcal{L}$$ over $L_a \times L_b \times L_c$, which determines in turn a map $$q_* f \to \pi^* \mathcal{L}.$$ Since a map of $(\Z/2)^3$-torsors is always an isomorphism, $$ q_* f \cong \pi^* \mathcal{L} .$$

\end{proof}

\tpoint{Corollary}\label{Rpoint_implies_TripleMassey0}{ \em  For all $a,b,c$ in $R^*$, if $X(a,b,c)(R) \neq \emptyset$, then $\langle \kappa(a),\kappa(b), \kappa(c)  \rangle $ contains $0$.}

\begin{proof}
The map $\Et(\Spec R) \to B \pi_1^{\et}(R)$ induces an isomorphism $$\rH^1(\pi_1^{\et}(\Spec R), \Z/2) \to \rH^1(\Spec R, \Z/2).$$ Thus it suffices (in fact it is equivalent) to show that $\langle \kappa(a),\kappa(b), \kappa(c)  \rangle$ contains $0$ in $\rH^2(\pi_1^{\et}(\Spec R), \Z/2)$.

Choose $\Spec R \to X(a,b,c)$, and let $E$ denote the pull back of the $U_4$ torsor $\G_{m, R}^4 \times \mathcal{L} \to \mathcal{X}$ of Theorem \ref{have_torsors} to $\Spec R$ via $$\Spec R \to X(a,b,c) \to \mathcal{X}.$$ By Theorem \ref{have_torsors}, $q_* E$ is isomorphic to the pullback of $\mathcal{L}$ along $$(a,b,c): \Spec R \to  \Spec R[a,c, b, a^{-1}, b^{-1},c^{-1}].$$ Viewing $E$ as an element of $\rH^1(\Spec R, U_4) \cong  \Hom (\pi_1^{\et}(\Spec R), U_4)$,  we have a homomorphism $$\pi_1^{\et}(\Spec R) \to U_4 $$ such that the composition $$\pi_1^{\et}(\Spec R) \to U_4 \to (\Z/2)^3$$ is $(\kappa(a), \kappa(b), \kappa(c))$.\hidden{$\kappa(a)$ determines a well-defined cocycle $\pi_1^{\et}(\Spec R) \to \Z/2$.} Thus  $\langle \kappa(a),\kappa(b), \kappa(c)  \rangle $ contains $0$ in $\rH^2(\pi_1^{\et}(\Spec R), \Z/2)$ by \ref{MPGroup_cohomology}.\end{proof}

\epoint{Example}\label{Gartner} The following example is work of Jochen G\"artner.  Let $T= \{l_1,l_2,\ldots, l_n,2,\infty \}$ be a set of odd primes union $\{2,\infty\}$, $\Q_T(2)$ denote the maximal $2$-extension of $\Q$ unramified outside $T$ and let $G_T(2) = \Gal(\Q_T(2)/\Q)$. For an appropriate choice of topological generators for $G_T(2)$ (similar to \cite[2.1.3 p 47 $\tau_i$]{Vogel_thesis} and \cite[3.2.5 $\tau_i$]{Gaertner_thesis}), the dual basis for $\rH^1(G_T(2), \Z/2)$ (\cite[p 49 $\chi_i$]{Vogel_thesis} \cite[3.3.5]{Gaertner_thesis}) contains the $\kappa(l_i)$. For $l_i$ such that the Legendre symbols satisfy $$\left( \frac{l_i}{l_j}\right) = 1\text{ for }i \neq j$$ $$l_i \equiv 1 \mod 8,$$ taking the cup product with $\kappa(l_i)$ produces the $0$ map $\rH^1(G_T(2), \Z/2) \to \rH^2(G_T(2), \Z/2)$. This gives a uniquely defined Massey product $$\langle -,-,- \rangle: \{ \kappa(l_1), \kappa(l_2), \ldots, \kappa(l_n)\}^3 \to \rH^2(G_T(2), \Z/2).$$ With results of Vogel and Morishita computing the triple Massey product in terms of R\'edei symbols \cite[Th. 2.1.16]{Vogel_thesis}, G\"artner can show using MAGMA that for $n=3$ and $(l_1,l_2,l_3) = (313,457,521)$, the Massey product $\langle \kappa(l_1), \kappa(l_2), \kappa(l_3) \rangle$ is not zero. It follows that there is no homomorphism $G_T(2) \to U_4$ such that composing with $a_{12} \times a_{23} \times a_{34}$ represents $(\kappa(l_1), \kappa(l_2), \kappa(l_3))$. Since $U_4$ is a $2$-group and $G_T(2)$ is the maximal $2$-quotient of $\pi_1(\Spec \Z[\frac{1}{l_1}, \frac{1}{l_2}, \frac{1}{l_3}, \frac{1}{2}])$, this implies that there is no such homomorphism $$\pi_1(\Spec \Z[\frac{1}{l_1}, \frac{1}{l_2}, \frac{1}{l_3}, \frac{1}{2}]) \to U_4.$$ It follows by Corollary \ref{Rpoint_implies_TripleMassey0} that $X(313,457,521)$ has no $\Z[\frac{1}{313}, \frac{1}{457}, \frac{1}{512}, \frac{1}{2}]$-points.

Now consider the case where $R$ is a field $F$.

\tpoint{Theorem}\label{SplittingVarietyThm}{\em Let $F$ be a field of characteristic $\neq 2$. For all $a,b,c$ in $F^*$, the scheme $X(a,b,c)$ has an $F$-point if and only if $\langle \kappa(a),\kappa(b), \kappa(c)  \rangle $ contains $0$.}

So $X(a,b,c)$ is a splitting variety for $\langle \kappa(a),\kappa(b), \kappa(c)  \rangle $. The following proof is similar to \cite[Prop. 4.11]{Berhuy_Favi}, and we thank Burt Totaro for sending us this reference.

\begin{proof}
Assume that $\langle \kappa(a),\kappa(b), \kappa(c)  \rangle$ contains $0$. We show that $X(a,b,c)$ has an $F$-point, and this is sufficient by Corollary \ref{Rpoint_implies_TripleMassey0}.

The triple Massey product $\langle \kappa(a),\kappa(b), \kappa(c)  \rangle$ in continuous group cohomology of $\pi_1^{\et}(\Spec F) \cong \Gal(F^s/ F)$ contains $0$, where $F^s$ denotes a separable closure of $F$. By Remark \ref{MPGroup_cohomology}, we have an element of $\rH^1(\Spec F, U_4)$ whose image under $q_* : \rH^1(\Spec F, U_4) \to \rH^1(\Spec F, \Z/2)^3$ is $(\kappa(a), \kappa(b), \kappa(c)).$ Let $$\sigma:\Gal(F^s/ F) \to U_4$$ be a homomorphism representing this element. 

Choose square roots $\alpha$, $\beta$, and $\gamma$ of $a$, $b$, and $c$ respectively in $F^s$. Then $$a_{12} \sigma (g) = (g \alpha) / \alpha \in \mu_2(F^s) \cong \Z/2$$ for all $g$ in $\Gal(F^s/F)$. Similarly $a_{23} \sigma (g) = (g \beta) / \beta $, and  $a_{34} \sigma (g) = (g \gamma) / \gamma $.

Let $\rho: U_4 \to \GL_4$ be the homomorphism corresponding to the representation $V = \operatorname{Ind}_{\tilde{H}}^{U_4} \sigma_{14}$ from the proof of Theorem \ref{have_torsors}. The image of $\sigma$ under $\rho_*: \rH^1(\Spec F, U_4 ) \to \rH^1(\Spec F, \GL_4)$ is trivial since $\rH^1(\Spec F, \GL_4 )$ is the pointed set with one element. Thus we have a matrix $A$ in $\GL_4 F^s$ such that $$\rho \sigma (g) = A^{-1} (gA) ,$$ for all $g$ in $\Gal(F^s/F)$.

The kernel of a non-zero $F$-linear map $F^4 \to F^s$ has dimension $< 4$. For $F$ an infinite field, the $F$-vector space $F^4$ is not contained in a union of finitely many dimension $<4$ sub-vector spaces. \hidden{claim: A vector space $V$ over an infinite field $F$ is not a finite union of positive codimension sub-vector spaces. proof:  By induction on the dimension of $V$. If $V$ is dimension $1$, then all positive codimension sub-vector spaces are $\{0\}$ and the claim follows. Assume $\dim V > 1$, and fix a finite collection of positive codimension sub-vector spaces of $V$. Each of these sub-vector spaces is contained in the kernel of an element of the dual $V^*$ of $V$. Let $I$ be a finite set of elements of $V^*$ such that $\cup_{i \in I} \Ker i$ contains the original finite collection of positive codimension sub-vector spaces. The projectivization $\proj V^*$ of the dual of $V$ has infinitely many elements. Thus there exists $\varphi$ in $\proj V^*$ not contained in the image of $I$ in $\proj V^*$. Thus the restriction to $\Ker \varphi$ of every element of $I$ is a non-zero functional. By induction, $$\cup_{i \in I} \Ker i\vert_{\Ker \varphi} \neq \Ker \varphi,$$ and it follows that $$\cup_{i \in I} \Ker I \neq V $$ as desired. }Thus there exists $\mu = (\mu^*_1, \mu^*_2, \mu^*_3, \mu^*_4)$ in $F^4$ such that $A^{-1} \mu$ is in $({F^s}^*)^4$. If $F$ is a finite field of characteristic $> 2$, there also exists such a $\mu$. To see this: 


If $F$ is a finite field of $q > 4$ elements, then the number of elements in the union of four positive codimension sub-vector spaces of $F^4$ is  less than $4 q^3< q^4$, whence there exists $\mu$. 


Otherwise $F$ has $q=3$ elements.

For $S \subset \{1,2,3,4\}$, let $u_S: (F^s)^4 \to (F^s)^{\vert S \vert}$ denote the projection onto the coordinate axes contained in $S$. 

First suppose that $\dim_F \Ker u_i A^{-1} < 3$ for some $i$. The number of elements of $\cup_{j \neq i} \Ker u_j A^{-1} $ is less than or equal to $q^3 + (q^3 - q^2) + (q^3 - q^2) $ because two dimension-$3$ sub-vector-spaces must intersect in a plane, and any positive codimension sub-vector-space can be enlarged to be dimension $3$. Thus the number of elements of $\cup_{j} \Ker u_j A^{-1} $ is less than or equal to $q^3 + (q^3 - q^2) + (q^3 - q^2) + q^2  = 3 q^3 -q^2 < q^4$, so there exists $\mu$.

We can thus suppose that  $\dim_F \Ker u_i A^{-1} = 3$ for all $i$. Since $\Ker u_S A^{-1} = \cap_{i \in S} \Ker u_i A^{-1}$, the dimension of $\Ker u_S A^{-1}$ is $\geq 4 - \vert S \vert$.

Take $N \leq 4$. An $F$-linear map $L:F^4 \to (F^s)^{N}$ determines an $F^s$-linear map $L \otimes F^s: (F^s)^4 \to (F^s)^{N}$, and $\dim_F L(F^4) \geq \dim_{F^s} (L \otimes F^s) ((F^s)^4)$. Thus $$ u_S A^{-1}: F^4 \to (F^s)^{\vert S \vert}$$ has rank $\geq \vert S \vert$ and kernel of dimension $\leq 4- \vert S \vert$.  Thus $\dim \Ker u_S A^{-1} = 4 - \vert S \vert$, and $$\vert \cup_{i = 1}^4 \Ker u_i  A^{-1} \vert = 4 q^{3} - {4 \choose 2} q^2 + {4 \choose 3} q - 1 < q^4.$$

Thus there exists $\mu = (\mu^*_1, \mu^*_2, \mu^*_3, \mu^*_4)$ in $F^4$ such that $A^{-1} \mu$ is in $({F^s}^*)^4$.

$A^{-1} \mu \times (\alpha, \beta, \gamma)$ determines an $F^s$-point of $\G_m^4 \times \mathcal{L}$. Furthermore, \begin{align*}g^{-1} (A^{-1} \mu \times (\alpha, \beta, \gamma)) & = \rho \sigma (g) A^{-1} \mu \times (g^{-1} \alpha, g^{-1} \beta, g^{-1} \gamma)  = \rho \sigma (g) A^{-1} \mu \times (\kappa(a) \alpha, \kappa(b) \beta, \kappa(c) \gamma)\\ & =  \rho \sigma (g) A^{-1} \mu \times q \sigma (g) (\alpha, \beta, \gamma) \end{align*} for all $g$ in $\Gal(F^s/F)$, so $A^{-1} \mu \times (\alpha, \beta, \gamma)$ determines an $F$-point of $X(a,b,c)$.

\end{proof}

\section{Vanishing for global fields}

Let $F$ be a field of characteristic $\neq 2$, and choose $a,c$ in $F^*$. Let $N: F^4 \to F$ be defined $$N(y_1, y_2, y_3, y_4) = (y_1^2- a y_2^2 + c y_3^2 - ac y_{4}^2 )^2 - c(2 y_1 y_{3} - 2 a y_{2} y_{4} )^2.$$

 Let $\mu$ be any prime of $F[\alpha, \gamma]/\langle \alpha^2 - a, \gamma^2 - c \rangle$, and let $F(\mu) \cong F[\sqrt{a}, \sqrt{c}]$ denote the residue field at $\mu$. Let $N_{F(\mu)/F} : F(\mu)^* \to F$ denote the norm map.

\tpoint{Proposition}\label{Norms_compatible}{\em The subsets $N(F^4 -\{0\} )$ and $N_{F(\mu)/F} (F(\mu)^*)$ of $F^*$ are equal.}

\begin{proof}
In $F(\mu)$, we have the equality $$N(y_1, y_2, y_3, y_4) = \prod_{i,j \in \{0,1\}} (y_1+  (-1)^{i} y_2 \sqrt{a} + (-1)^j y_3 \sqrt{c} + (-1)^{i+j} y_4 \sqrt{a}\sqrt{c} ).$$ If the extension $F \subseteq F(\mu)$ has degree $4$, then $N (y_1, y_2, y_3, y_4)= N_{F(\mu)/F}(y_1 + y_2 \sqrt{a}+ y_3 \sqrt{c} + y_4\sqrt{a}\sqrt{c})$, giving the result immediately. 

If $F \subseteq F(\mu)$ has degree $2$, then we may assume $F(\mu) \cong F[\sqrt{a}]$ for $a \neq 1$ in $F^*/{F^*}^2$, so $N(y_1, y_2, y_3, y_4) = $ $$ N_{F(\mu)/F}(y_1 + y_2 \sqrt{a}+ y_3 \sqrt{c} + y_4\sqrt{a}\sqrt{c})N_{F(\mu)/F}(y_1 + y_2 \sqrt{a}- y_3 \sqrt{c} - y_4\sqrt{a}\sqrt{c}).$$ Since the map $F^4 \to F(\mu)^2 $ given by $$(y_1, y_2, y_3, y_4) \mapsto (y_1 + y_2 \sqrt{a}+ y_3 \sqrt{c} + y_4\sqrt{a}\sqrt{c}, y_1 + y_2 \sqrt{a}- y_3 \sqrt{c} - y_4\sqrt{a}\sqrt{c}) $$ is surjective and $N_{F(\mu)/F} (1) = 1$ , we have that $N(F^4 -\{0\} ) \supseteq  N_{F(\mu)/F} (F(\mu)^*)$. Since  $N_{F(\mu)/F}$ is multiplicative, $N(F^4 -\{0\} ) \subseteq  N_{F(\mu)/F} (F(\mu)^*)$, giving the result.

If $F \subseteq F(\mu)$ has degree $1$, then $N_{F(\mu)/F}$ is the identity map and $N_{F(\mu)/F} (F(\mu)^*)= F^*$. Since the linear map $F^4 \to F^4$ determined by \begin{equation*}\left( \begin{array}{cccc}
 1&   \sqrt{a} &  \sqrt{c} & \sqrt{a}\sqrt{c} \\
 1&  \sqrt{a} & - \sqrt{c} & - \sqrt{a}\sqrt{c} \\
 1&  - \sqrt{a} &  \sqrt{c} & - \sqrt{a}\sqrt{c} \\
  1&  - \sqrt{a} & - \sqrt{c} &  \sqrt{a}\sqrt{c} \end{array} \right) \end{equation*} is an isomorphism, $N(F^4 -\{0\} ) = F^*$ as well.
\end{proof}

From the definition of $X(a,b,c)$ and Proposition \ref{Norms_compatible} we obtain:

\tpoint{Corollary}\label{Nimage=Normimage}{\em $X(a,b,c)(F) \neq \emptyset$ if and only if for any prime $\mu$ of  $F[\alpha, \gamma]/\langle \alpha^2 - a, \gamma^2 - c \rangle$, there exists $x \in F^*$ and $y \in F(\mu)^*$ such that $$b x^2 = N_{F(\mu)/F}(y).$$}

The following corollary can also be shown directly, using for instance \cite[Lemma 16]{WPIA} \hidden{$PIAc$ such that $D PIAc= MPSV b \cup a $, $PIAb = MPSVa$}, but we include a proof using the splitting variety $X(a,b,a)$.

\tpoint{Corollary}\label{<a,b,a>=0}{\em Let $F$ be a field of characteristic $\neq 2$, and $a,b \in F^*$. The triple Massey product $\langle \kappa(a),\kappa(b), \kappa (a) \rangle $ contains $0$ whenever it is defined.}

\begin{proof}
Let $L = F[\sqrt{a}]$. When $\langle \kappa(a),\kappa(b), \kappa (a) \rangle $ is defined, $a \cup b = 0$ which implies that $b$ is trivial in $F^*/((N_{L/F} F^*){F^*}^2)$, for instance because $a x^2 + b y^2 = z^2$ is a splitting variety for the cup-product. By Corollary \ref{Nimage=Normimage}, this implies that $X(a,b,a)\neq \emptyset$. Thus $\langle \kappa(a),\kappa(b), \kappa (a) \rangle $ contains $0$ by Theorem \ref{SplittingVarietyThm}.
\end{proof}

The splitting variety $X(a,b,c)$ satisfies the Hasse principle.

\tpoint{Theorem}\label{HasseX}{\em Let $F$ be a global field of characteristic $\neq 2$ and $a,b,c \in F^*$. If $X(a,b,c)(F_{\nu}) \neq \emptyset$ for all places $\nu$ of $F$, then $X(a,b,c)(F) \neq \emptyset$.}

\begin{proof}
Let $L = F[\sqrt{a}, \sqrt{c}]$.

Choose a place $\nu$ of $F$. By Corollary \ref{Nimage=Normimage}, we have that $b$ is in $F_{\nu}^2 N_{L_{\mu}/F_{\nu}}(L_\mu^*) $ for every place $\mu$ of $L$ above $\nu$.

By the Hasse norm theorem mod squares of Leep and Wadsworth \cite[Thm 4.5]{LW_transfer}, $$\{ b \in F^*: b \in {F_{\nu}^*}^2 N_{L_{\mu}/F_{\nu}}(L_{\mu}^*)~~ \forall \mu \vert \nu\} \subset \{  b \in F^*: b \in {F^*}^2 N_{L/F}(L^*)\},$$ showing that $X(a,b,c)(F) \neq \emptyset$ by Corollary \ref{Nimage=Normimage}.
\end{proof}

\tpoint{Lemma}\label{localM=0}{\em Let $F$ be a local field of characteristic $\neq 2$ and $a,b,c \in F^*$. The triple Massey product $\langle \kappa(a),\kappa(b), \kappa(c) \rangle$ contains $0$ whenever it is defined.}

\begin{proof}

It is sufficient to show that $$ \kappa(a) \cup (-) : \rH^1(\Spec F, \Z/2) \to  \rH^2(\Spec F, \Z/2) $$  is surjective by \ref{indeterminacy_remark}. The $2$-torsion $\Br(F)[2]$ of the Brauer group of $F$ is isomorphic to $\rH^2(\Spec F, \mu_2)$ by Hilbert 90 and the short exact sequence $$1 \to \mu_2 \to \G_m \stackrel{z \mapsto z^2}{\to} \G_m \to 1.$$ Since $F$ is a local field, $\Br (F) \cong \Q/\Z$, whence $\rH^2 (\Spec F, \mu_2) \cong \Z/2$ (see \cite[VI]{CF}, for instance), so $\rH^2(\Spec F, \Z/2) \cong \Z/2$. It thus suffices to see that the cup-product pairing $$  \rH^1(\Spec F, \Z/2) \otimes  \rH^1(\Spec F, \Z/2) \to  \rH^2(\Spec F, \Z/2) \cong \Z/2$$ is non-degenerate. This is true by Tate duality \cite[Thm 7.2.6]{coh_num_fields} \cite{Tate_ICM}.

\end{proof}

\tpoint{Theorem}\label{globalM=0}{\em Let $F$ be a global field of characteristic $\neq 2$ and $a,b,c \in F^*$. The triple Massey product $\langle \kappa(a),\kappa(b), \kappa(c) \rangle$ contains $0$ whenever it is defined.}

\begin{proof}
Suppose $\kappa(a) \cup \kappa(b) = \kappa(b) \cup \kappa(c) = 0$, so that $\langle \kappa(a),\kappa(b), \kappa(c) \rangle$ is defined. By Lemma \ref{localM=0} and Theorem \ref{SplittingVarietyThm}, we have $X(a,b,c) (F_{\nu}) \neq \emptyset$. By Theorem \ref{HasseX}, we have $X(a,b,c)(F) \neq \emptyset$, so $\langle \kappa(a),\kappa(b), \kappa(c) \rangle$ contains $0$ by another application of Theorem \ref{SplittingVarietyThm}.
\end{proof}

\epoint{Remark}\label{local_global_Galois_remark} By the Hasse--Brauer--Noether local-global principle $$0 \to \rH^2(\Spec F, \Z/2) \to \bigoplus_{\nu} \rH^2(\Spec F_{\nu},\Z/2\Z)  \xrightarrow{\sum_v \operatorname{inv}_v} \frac{1}{2}\Z/\Z \to 0,$$ $\rH^2(-, \Z/2)$ of a global field injects into the direct sum of $\rH^2(-, \Z/2)$ of its places. This shows the local-global principle for cup products and thus the Hasse principle for their splitting varieties. 

The situation for higher Massey products is more subtle. Suppose $\langle \kappa(a), \kappa(b), \kappa(c)\rangle$ vanshes in $ \mathcal{C}^*(\Spec F_{\nu},\Z/2\Z)$ for all places $\nu$, where $ \mathcal{C}^*(\Spec F_{\nu},\Z/2\Z)$ denotes the differential graded algebra of continuous Galois cochains. This means that for all $\nu$, there exist $E^{\nu}_{ab}$ and $E^{\nu}_{bc}$ in $ \mathcal{C}^1(\Spec F_{\nu},\Z/2\Z)$ such that $\delta E^{\nu}_{ab} = \kappa(a) \cup\kappa( b)$, $\delta E^{\nu}_{bc} = \kappa(b) \cup \kappa(c)$, and $E^{\nu}_{ab} \cup \kappa(c) + \kappa(a) \cup E^{\nu}_{bc}$ vanishes in $\rH^2(\Spec F_{\nu},\Z/2\Z)$. If the $E^{\nu}_{ab}$ and $E^{\nu}_{bc}$ can be chosen compatibly, we could conclude the vanishing of $\langle \kappa(a), \kappa(b), \kappa(c)\rangle$ for $F$ by the local-global principle for $\rH^2(-,\Z/2)$. For instance, if the map \begin{equation}\label{FU3specialization} \rH^1(\Spec F, U_3) \to \oplus_{\nu} \rH^1(\Spec F_{\nu}, U_3)\end{equation} is surjective, we could choose the $E^{\nu}_{ab}$ and $E^{\nu}_{bc}$ compatibly and conclude vanishing. A $2$-nilpotent version of Poitou-Tate duality (see \cite[Ch 18]{Stix} for a non-abelian Poitou-Tate duality result) would give information about \eqref{FU3specialization}, however we do not know of an argument along these lines which shows that $\langle \kappa(a), \kappa(b), \kappa(c)\rangle$ contains $0$ in $\rH^2(\Spec F,\Z/2\Z)$ when it is defined and contains $0$ at all places. 

Instead of directly using Galois cohomology to study local-global properties for triple Massey products, Theorem \ref{HasseX} used the computed equation for $X(a,b,c)$.

} 

\bibliographystyle{MPSV}


\bibliography{MPSV}


\end{document}